\documentclass[11pt]{article}  

\usepackage[dvips]{epsfig}
\usepackage{latexsym}
\usepackage{amsfonts,amsmath,amssymb}
\newtheorem{corollary}{Corollary}

\newtheorem{example}{Example}
\newtheorem{theorem}{Theorem}
\newtheorem{proposition}{Proposition}

\newtheorem{definition}{Definition}
\newcommand{\qed}{\mbox{$\Diamond$}\vspace{\baselineskip}}
\newenvironment{proof}{\noindent{\bf Proof:}}{\qed}
\newcommand{\Var}{\hbox {Var}}

\usepackage{lineno}

\begin{document}
\title{The copies of any permutation pattern are asymptotically normal} 

\author{Mikl\'os B\'ona \\
Department of Mathematics \\
University of Florida\\
 Gainesville FL 32611-8105\\
bona@math.ufl.edu \thanks{Partially supported
by an NSA Young Investigator Award.}}

\date{}
\maketitle

\begin{abstract} We prove that the number of copies of any given permutation
pattern $q$ has an asymptotically normal distribution in random permutations. 
\end{abstract}

\section{Introduction}
The classic definition of pattern avoidance for permutations is as follows.
Let $p=p_1p_2\cdots p_n$ be a permutation, let $k<n$, and let $q=q_1q_2\cdots
q_k$ be another permutation. We say that $p$ {\em contains} $q$ as a pattern
if there exists a subsequence $1\leq i_1<i_2<\cdots <i_k\leq n$
 so that for all indices
$j$ and $r$, the inequality $q_j<q_r$ holds if and only if the inequality
$p_{i_j}<p_{i_r}$ holds. If $p$ does not contain $q$, then we say
that $p$ {\em avoids} $q$.  In other words, $p$ contains $q$ if $p$ has
a subsequence of entries, not necessarily in consecutive positions, which
relate to each other the same way as the entries of $q$ do.

In a recent survey paper \cite{monsurvey} on the monotone permutation
pattern $12\cdots k$, we have shown that if $X_n$ is the random variable
counting copies of that pattern in a randomly selected permutation
of length $n$, then as $n$ goes to infinity, $X_n$ converges (in distribution)
to a normal distribution. When we say ``random permutation'', we mean that
each permutation of length $n$ is selected with probability $1/n!$.

In this paper, we will generalize that result for {\em any} permutation
pattern $q$, and the variable $X_{n,{\bf q}}$ counting the copies of $q$ in
permutations of length $n$. The proof is very similar to the monotone case; 
just some details
have to be modified. The result is a far-reaching generalization of
the classic results (see \cite{fulman}) for more references) that descents
and inversions of random permutations are asymptotically normal.
 As a byproduct, we will see how close $\Var(X_{n,{\bf q}})$
and $\Var(X_{n,{\bf 12\cdots k}})$
 are to each other, for any pattern $q$ of length
$k$. 

\section{The Proof of Our Theorem}
\subsection{Background and Definitions}
We need to introduce some notation for transforms of the random variable
$Z$. Let $\bar{Z}=Z-E(Z)$, let $\tilde{Z}=\bar{Z}/\sqrt{\Var( Z)}$, and let
$Z_n\rightarrow N(0,1)$ mean that $Z_n$ converges in distribution to the 
standard normal variable. 

\begin{definition}
Let $\{Y_{n,k}|k=1,2,\cdots ,N_n\}$ be an array of
 random variables.
 We say that a graph $G$ is 
a {\em dependency graph} for   $\{Y_{n,k}|k=1,2\cdots , N_n\}$
 if the following
two conditions are satisfied:
\begin{enumerate}
\item There exists a bijection between the random variables $Y_{n,k}$ and
the vertices of $G$, and
\item If $V_1$ and $V_2$ are two disjoint sets of vertices of $G$ so that
no edge of $G$ has one endpoint in $V_1$ and another one in $V_2$, then
the corresponding sets of random variables are independent.
\end{enumerate}
\end{definition}

Note that  the dependency graph of a
 family of variables is not unique. Indeed if $G$ is a dependency graph
for a family and $G$ is not a complete graph,
 then we can get other dependency graphs for the family
by simply adding new edges to $G$. 

Now we are in position to state Janson's theorem, the famous
{\em Janson dependency criterion}.

\begin{theorem} \cite{janson} \label{janson}
Let $Y_{n,k}$ be an array of random variables such that for all $n$, and
for all $k=1,2,\cdots ,N_n$, the inequality $|Y_{n,k}|\leq A_n$ holds for
some real number $A_n$, and that the maximum degree of a dependency
graph of $\{Y_{n,k} | k=1,2,\cdots ,N_n \}$ is $\Delta_n$. 

Set $Y_n=\sum_{k=1}^{N_n} Y_{n,k}$ and $\sigma_n^2= \Var ( Y_n)$. If there
is a natural number $m$ so that
\begin{equation} \label{jansencond}
N_n\Delta_n^{m-1} \left (\frac{A_n}{\sigma_n} \right )^m \rightarrow 0,
\end{equation}
as $n$ goes to infinity, then \[ \tilde{Y}_n \rightarrow N(0,1) .\]
\end{theorem}

\subsection{Verifying the Conditions of Janson's Criterion}
Let $q$ be a fixed pattern of length $k$. As $q$ is fixed for the rest
of this paper, we will  mark our variables $X_n$ instead of $X_{n,{\bf q}}$,
in order to avoid excessive indexing. 
  
Let us order the ${n\choose k}$ subwords of length $k$ of the permutation
$p_1p_2\cdots p_n$ linearly in some way. 
For $1\leq i\leq {n\choose k}$, let $X_{n,i}$ 
 be the indicator random
variable of the event that in a randomly selected permutation of length $n$,
the $i$th subword of length $k$ in the permutation $p=p_1p_2\cdots p_n$
is a $q$-pattern. We will now verify that the family of the
$X_{n,i}$ satisfies all conditions of the Janson Dependency Criterion.

First, $|X_{n,i}|\leq 1$ for all $i$ and all $n$, since the $X_{n,i}$ are 
indicator random variables. So we can set $A_n=1$. Second, $N_n={n\choose k}$,
the total number of subwords of length $k$ in $p$. Third, if $a\neq b$, then
$X_a$ and $X_b$ are independent unless the corresponding subwords intersect.
For that, the $b$th subword must intersect the $a$th subword in $j$ entries, 
for some $1\leq j\leq k-1$. For a fixed $a$th subword, the number of 
ways that can happen is $\sum_{j=1}^{k-1} {k\choose j}{n-k\choose k-j}=
{n\choose k}-{n-k \choose k}-1$, where we used 
 the well-known Vandermonde identity to compute the sum.
Therefore, 
\begin{equation} \label{maxdegree} 
\Delta_n \leq {n\choose k}-{n-k \choose k}-1.
\end{equation}
In particular, note that (\ref{maxdegree}) provides an upper bound for
$\Delta_n$ in terms of a polynomial function of $n$ that is 
of degree $k-1$ since terms of degree
$k$ will cancel.

There remains the task of finding a lower bound for $\sigma_n$ that 
we can then use in applying Theorem \ref{janson}. Let $X_n=
\sum_{i=1}^{n\choose k} X_{n,i}$. We will show the following.

\begin{proposition} \label{varprop}
 There exists a positive constant $c$ so that
for all $n$, the inequality
\[\Var(X_n)\geq cn^{2k-1}\]
holds.
\end{proposition}

\begin{proof}
By linearity of expectation, we have
\begin{eqnarray} \label{variance}
\Var (X_n) & = & E(X_n^2) - (E(X_n))^2 \\
 & = & E \left (\left( \sum_{i=1}^{{n\choose k}} X_{n,i} \right )^2 \right )
-  \left (E \left (\sum_{i=1}^{{n\choose k}} X_{n,i} \right ) \right )^2 \\
 & = & E \left (\left( \sum_{i=1}^{{n\choose k}} X_{n,i} \right )^2 \right )
- \left( \sum_{i=1}^{{n\choose k}} E(X_{n,i}) \right )^2  \\
 \label{lastone} & = &  \sum_{i_1, i_2}
E(X_{n,i_1}X_{n,i_2})  - \sum_{i_1, i_2}
E(X_{n,i_1})E(X_{n,i_2}).
\end{eqnarray}

Let $I_1$ (resp. $I_2$) denote the $k$-element subword of $p$ indexed
by $i_1$, (resp. $i_2$). Clearly, it suffices to show that
\begin{equation} \label{simplified} \sum_{|I_1\cap I_2| \leq 1}
E(X_{n,i_1}X_{n,i_2}) - \sum_{i_1, i_2}
E(X_{n,i_1})E(X_{n,i_2}) \geq cn^{2k-1},\end{equation}
since the left-hand side of (\ref{simplified}) is obtained from the 
(\ref{lastone}) by removing the sum of some positive terms, that is,
the sum of all $E(X_{n,i_1}X_{n,i_2})$ where  $|I_1\cap I_2| >1$.

As $E(X_{n,i})=1/k!$ for each $i$, the sum with negative sign in 
(\ref{lastone}) is
\[ \sum_{i_1, i_2}
E(X_{n,i_1})E(X_{n,i_2}) ={n\choose k}^2 \cdot \frac{1}{k!^2},\]
which is a polynomial function 
 in $n$, of degree $2k$ and of leading coefficient
$\frac{1}{k!^4}$. As far as the summands in  (\ref{lastone}) with a positive
sign go, {\em most} of them are also equal to  $\frac{1}{k!^2}$. More
precisely, $E(X_{n,i_1}X_{n,i_2})=\frac{1}{k!^2}$ when 
 $I_1$ and $I_2$ are disjoint, and that happens for 
${n\choose k}{n-k\choose k}$ ordered pairs $(i_1,i_2)$
 of indices. The sum of these
summands is
\begin{equation}
\label{disjoint} d_n={n\choose k}{n-k\choose k} \frac{1}{k!^2},
\end{equation}
which is again a polynomial function in $n$, of degree $2k$ and with leading 
coefficient
$\frac{1}{k!^4}$. So summands
 of degree $2k$ will cancel out in (\ref{lastone}). (We will see in the next
paragraph that the summands we have not yet considered add up to a polynomial
of degree $2k-1$.)

In fact, considering the two types of summands we studied in
(\ref{lastone}) and (\ref{disjoint}), we see that they add up to 
\begin{eqnarray} 
{n\choose k}{n-k\choose k} \frac{1}{k!^2}-{n\choose k}^2  \frac{1}{k!^2}
& = & n^{2k-1} \frac{2{k\choose 2}-{2k-1\choose 2}}{k!^4}+O(n^{2k-2}) \\
\label{theeasy} & = & n^{2k-1} \frac{-k^2}{k!^4} +O(n^{2k-2}) .
\end{eqnarray}

Next we look at ordered pairs of indices $(i_1,i_2)$ so that the corresponding
subwords $I_1$ and $I_2$ intersect in exactly one entry, the entry
 $x$. Let us restrict our attention to the special case when
  $I_1$ and $I_2$ both form $q$-patterns, and 
$x$ is the $a$th smallest entry in $I_1$ and the $b$th smallest entry in
 $I_2$. Given $q$, the pair $(a,b)$ 
describes the location of $x$ in $I_1$ and in $I_2$
as well. Let $I_1'$ (resp. $I_2'$) denote the set of
 $a-1$ positions in $I_1$ (resp. 
$b-1$ positions in 
$I_2$) which must contain entries smaller than $x$ given that
 $I_1$ (resp. $I_2$)  forms a $q$-pattern. Similarly, let  
$I_1''$ (resp. $I_2''$) denote the
set of $k-a$ positions in $I_1$ (resp. 
$k-b$ positions $I_2$) which must contain entries larger than $x$ given that
 $I_1$ (resp. $I_2$) 
 forms a $q$-pattern.

\begin{example} 
Let $q=35142$, and let us say that $I_1$ and $I_2$ both form $q$-patterns,
and they intersect in one entry $x$ that is the third smallest entry in 
$I_1$ and the fourth smallest entry in $I_2$ (so $a=3$, and $b=4$).
Then $x$ is the leftmost entry of $I_1$ and the next-to-last entry of $I_2$.
Furthermore, the third and fifth positions of $I_1$ form $I_1'$ and the
second and fourth positions of $I_1$ form $I_1''$. Similarly, 
the first, third, and fifth positions of $I_2$ form $I_2'$ and the second
position of $I_2$ forms $I_2''$.  
\end{example}

Let $q_a$ (resp. $q_b$) be the pattern obtained from $q$ by removing its
$a$th smallest (resp. $b$th smallest) entry. 

Note that  $X_{i_1}X_{i_2}=1$ if and only if all of the following independent 
events hold.
\begin{enumerate}
\item In the $(2k-1)$-element set of entries that belong to $I_1\cup I_2$, 
the entry $x$ is the $(a+b-1)$th smallest. This happens with
probability $1/(2k-1)$. 
\item The $a+b-2$ entries in positions belonging to $I_1'\cup I_2'$
 must all be 
smaller than the $2k-a-b$ entries in positions belonging to 
 $I_1''\cup I_2''$. This happens with probability 
$\frac{1}{{2k-2\choose a+b-2}}$.
\item 
\begin{itemize}
\item the subword $I_1'$ is a pattern that is isomorphic to the
pattern formed by the $a-1$ smallest entries of $q$, 
\item the subword
$I_2'$ is a pattern that is isomorphic to the
pattern formed by the $b-1$ smallest entries of $q$, 
\item the subword $I_1''$ is a pattern that is isomorphic to the
pattern formed by the $k-a$ largest entries of $q$, and 
\item 
the subword $I_2''$ is a pattern that is isomorphic to the
pattern formed by the $k-b$ largest entries of $q$.
This happens with probability
$\frac{1}{(a-1)!(b-1)!(k-a)!(k-b)!}$.
\end{itemize}
\end{enumerate}

Therefore, if  $|I_1\cap I_2|=1$, then
\begin{eqnarray} \label{oneprob}
P(X_{i_1}X_{i_2}=1) & = & 
\frac{1}{(2k-1){2k-2\choose a+b-2}(a-1)!(b-1)!(k-a)!(k-b)!} \\
& = & \frac{1}{(2k-1)!}\cdot{a+b-2 \choose a-1}{2k-a-b\choose k-a}
.\end{eqnarray}

How many such ordered pairs $(I_1,I_2)$ are there? There are ${n\choose 2k-1}$
choices for the underlying set $I_1\cup I_2$.
Once that choice is made,
the $a+b-1$st smallest entry of $I_1\cup I_2$ will be $x$.
Then
the number of choices for the set of  entries other than $x$
that will be part of $I_1$ is ${a+b-2\choose a-1}{2k-a-b\choose k-a}$.
 Therefore, summing over all $a$ and $b$ and
recalling (\ref{oneprob}), 
\begin{eqnarray} \label{contribution}
p_n & = &  \sum_{|I_1\cap I_2|=1}P(X_{i_1}X_{i_2}=1)=\sum_{|I_1\cap I_2|=1}
E (X_{i_1}X_{i_2})  \\
\label{thehard} &  = &  
 \frac{1}{(2k-1)!}{n\choose 2k-1}\sum_{1\leq a,b\leq k}{a+b-2 \choose a-1}^2
{2k-a-b\choose k-a}^2.
\end{eqnarray}

The expression we just obtained is a polynomial of degree $2k-1$, in the
variable $n$. We claim that its leading coefficient is 
larger than $k^2/k!^4$. If we can show that, the proposition will be proved
since (\ref{theeasy}) shows that the summands not included in 
(\ref{contribution}) contribute about $-\frac{k^2}{k!^4}n^{2k-1}$ to 
the left-hand side of (\ref{simplified}). 

Recall that by the Cauchy-Schwarz inequality, if $t_1,t_2,\cdots, t_m$
are non-negative real numbers, then 
\begin{equation}\label{schwarz}
\frac{\left(\sum_{i=1}^m t_i\right)^2}{m} \leq \sum_{i=1}^m t_i^2,
\end{equation}
where equality holds if and only if all the $t_i$ are equal.

Let us apply this inequality with the numbers ${a+b-2 \choose a-1}^2
{2k-a-b\choose k-a}^2$ playing the role of the $t_i$, where $a$ and $b$ 
range from 1 to $k$.
We get that 
\begin{equation}
\label{cauchy} \sum_{1\leq a,b \leq k}{a+b-2 \choose a-1}^2
{2k-a-b\choose k-a}^2  > \frac{\left (
\sum_{1\leq a,b\leq k} {a+b-2 \choose a-1}
{2k-a-b\choose k-a} \right)^2}{k^2}. \end{equation}
We will use Vandermonde's identity to compute the right-hand side. To that
end, we first compute the sum of summands with a {\em fixed} $h=a+b$.
We obtain
\begin{eqnarray}
\sum_{1\leq a,b\leq k} {a+b-2 \choose a-1}
{2k-a-b\choose k-a} &  = & \sum_{h=2}^{2k} \sum_{a=1}^k 
{h-2\choose a-1}{2k-h\choose k-a} \\
 &  = & \sum_{h=2}^{2k} {2k-2\choose k-1} \\
 &  = & (2k-1) \cdot {2k-2\choose k-1}.
\end{eqnarray}

Substituting the last expression into the right-hand side of (\ref{cauchy})
yields
\begin{equation} \label{estimate} \sum_{1\leq a,b \leq k}{a+b-2 \choose a-1}^2
{2k-a-b\choose k-a}^2 >  \frac{1}{k^2} \cdot (2k-1)^2 \cdot
{2k-2\choose k-1}^2.\end{equation}

\noindent Therefore, (\ref{contribution}) and (\ref{estimate}) imply that
\[p_n>\frac{1}{(2k-1)!}{n\choose 2k-1}\frac{(2k-1)^2}{k^2} 
{2k-2\choose k-1}^2.\]
As we pointed out after (\ref{contribution}), $p_n$ is a polynomial of 
degree $2k-1$ in the variable $n$. The last displayed inequality shows that
its leading coefficient is larger than
\[ \frac{1}{(2k-1)!^2} \cdot \frac{1}{k^2} \cdot \frac{(2k-2)!^2}{(k-1)!^4}
=\frac{k^2}{k!^4}  \] as claimed. 

Comparing this with (\ref{theeasy})
 completes the proof of our Proposition.
\end{proof}

We can now return to the application of Theorem \ref{janson} to our 
variables $X_{n,i}$. By Proposition \ref{varprop}, there is an absolute
constant $C$ so that $\sigma_n>Cn^{k-0.5}$ for all $n$. 
So (\ref{jansencond}) will be satisfied if we show that
there exists a positive integer $m$ so that
\[{n\choose k} (dn^{k-1})^{m-1} \cdot (n^{-k+0.5})^m<
dn^{-0.5m} 
 \rightarrow 0.\]
Clearly, any positive integer $m$ is a good choice. So we have proved the
following theorem.

\begin{theorem} \label{main}
 Let $q$ be a fixed permutation pattern of length $k$, 
and let $X_n$ be the random variable counting
 occurrences of $q$ in permutations of length $n$. 
Then $\tilde{X}_n\rightarrow N(0,1)$. In other words, $X_n$ is asymptotically
normal.
\end{theorem}

The following Corollary shows how close the variances of the numbers of copies
of two given patterns are to each other. 

\begin{corollary} For any pattern $q$ of length $k$, we have
\[Var(X_{n,{\bf q}})=c_kn^{2k-1}+O(n^{2k-2}),\]
where \[c_k=\frac{1}{(2k-1)!^2}\sum_{1\leq a,b \leq k}{a+b-2 \choose a-1}^2
{2k-a-b\choose k-a}^2 -\frac{k^2}{k!^4}  .\]
\end{corollary}

We point out that this does {\em not} mean that $Var(X_{n,{\bf q}})$ does not
depend on $q$. It does, and it is easy to verify that $\Var(X_{4,{\bf 123}})
\neq
\Var(X_{4,{\bf 132}})$. However, it is only the terms of degree at most $2k-2$
of $\Var(X_{n,{\bf q}})$ that depend on $q$. 
 
\begin{proof}
Note that in the proof of Theorem \ref{main}, we have not used anything
about the pattern $q$ apart from  its length. Our claim then follows from
comparing (\ref{theeasy}) and (\ref{thehard}).
\end{proof}

\end{document}